\documentclass[10 pt, leqno]{amsart}

\usepackage {amsfonts}
\usepackage{amsthm}
\usepackage{amssymb}
\usepackage{latexsym}
\usepackage{amsmath}
\usepackage{mathrsfs}
\usepackage{comment}
\usepackage{pifont}

\theoremstyle{plain}
\newtheorem{tw}{Theorem}[section]

\newtheorem {lem} [tw]{Lemma}
\newtheorem {prop}[tw] {Proposition}

\newtheorem{cor}[tw]{Corollary}

\theoremstyle{definition}
\newtheorem {deft}[tw] {Definition}

\newcommand{\Pois}{R_{\Vcont}}

\newcommand{\TPois}{R_{\repT}}
\newcommand{\sTPois}{R_{s,\repT}}

\newcommand{\bc} {\mathbb C}

\newcommand{\bn}{\mathbb N}

\newcommand{\be}{\mathbb E}

\newcommand{\bz}{\Bbb Z}

\newcommand{\alg} {\mathsf{A}}

\newcommand {\Pop}{\Delta_s(\Vcont)}
\newcommand {\PoT}{\Delta_s(\repT)}
\newcommand {\WoT}{W_s(\repT)}
\newcommand {\GaT}{\Gamma_s(\repT)}
\newcommand {\Ran} {{\textrm{Ran}}}

\newcommand {\id} {{\textup{id}}}
\newcommand {\Ker} {{\textup{Ker}}}

\newcommand{\Hil}{\mathsf{H}}

\newcommand{\Kil}{\mathsf{K}}

\newcommand{\Vcont}{\mathcal{V}}
\newcommand{\Fock}{\mathcal{F}}

\newcommand{\Wcont}{\mathcal{W}}

\newcommand{\La}{\Lambda}
\newcommand{\Nr}{\bn_0^r}

\newcommand{\wT}{\wt{T}}

\newcommand{\Toep}{\mathcal{T}_{\be}}
\newcommand{\repT}{\overrightarrow{T}}
\newcommand{\repV}{\overrightarrow{V}}

\newcommand{\Lin}{\textrm{Lin}}

\newcommand{\oner}{\{1,\ldots,r\}}

\newenvironment{rlist}
{

\begin{enumerate}}
{\end{enumerate}}

\newcommand{\la}{\lambda}
\newcommand{\ra}{\rangle}
\newcommand{\rra}{\right\rangle}
\newcommand{\lla}{\left\langle}

\newcommand{\Zr}{\bz^r}

\newcommand{\ot}{\otimes}

\newcommand{\ol}{\overline}

\newcommand{\wt}{\widetilde}

\numberwithin{equation}{section}


\begin{document}
\author{Adam Skalski}
\footnote{\emph{Permanent address of the author:} Department of Mathematics, University of \L\'{o}d\'{z}, ul. Banacha 22, 90-238
\L\'{o}d\'{z}, Poland.}
\address{Department of Mathematics and Statistics, Lancaster University, Lancaster LA1 4YF, United Kingdom }
\email{a.skalski@lancaster.ac.uk}

\title{\bf  On isometric dilations of product systems of $C^*$-correspondences and applications to
families of contractions associated to higher-rank graphs}
\begin{abstract}
\noindent Let $\be$ be a product system of $C^*$-correspondences over $\Nr$. Some sufficient conditions for the existence of a
not necessarily regular isometric dilation of a completely contractive representation of $\be$ are established and difference
between regular and $^*$-regular dilations discussed. It is in particular shown that a minimal isometric dilation is $^*$-regular
if and only if it is doubly commuting.
 The case of product systems associated with higher-rank
graphs is analysed in detail.
\end{abstract}

\maketitle

Classical multi-dimensional dilation theory (\cite{Nagy}) for Hilbert space operators is concerned with dilating tuples of
contractions to tuples of isometries or unitaries, preserving as many properties of the original family as possible. In
particular if the tuple with which we start consists of mutually commuting operators, it is desirable to obtain a commuting
dilation. Celebrated examples of S.\,Parrott, N.\,Varopoulos and others show that a joint dilation of three or more commuting
contractions to commuting isometries need not exist. In general it is difficult to decide whether a given commuting tuple has a
commuting isometric dilation. On the other hand the existence of so-called regular or $^*$-regular dilations (i.e.\ dilations
satisfying additional conditions with respect to products of the original contractions and their adjoints, see for example
\cite{Timot}) can be detected via simple conditions corresponding to positive-definiteness of certain operator-valued functions
associated with the initial tuple.

 In a recent  paper  \cite{graphdil} together with J.\,Zacharias we considered
dilations of $\La$-contractions, that is tuples of operators satisfying commutation relations encoded by a (higher rank) graph
$\La$.
 It has now become clear that using the constructions provided by I.\,Raeburn and A.\,Sims in \cite{toep} some of the
results of \cite{graphdil}  can be viewed as statements on completely contractive representations of the canonical product system
of $C^*$-correspondences associated to $\La$. Product systems of $C^*$-correspondences were first defined in \cite{fowl} as
generalisations of product systems of Hilbert spaces and quickly proved to provide a natural framework for extensions of the
classical multi-dimensional dilation theory to more complicated objects (see \cite{Solk} and references therein). The questions
about the existence of a joint dilation of a family of contractions satisfying certain commutativity relations to an analogous
family of isometries translates here into a question on the existence of an isometric dilation of a (completely) contractive
representation of a given product system.

Motivated by the observations above we show in this paper that the generalised Poisson transform constructed in \cite{graphdil}
(see also \cite{PPois}, \cite{MSker}) can be associated to a completely contractive representation of a  product system of
$C^*$-correspondences over $\Nr$, if only the system enjoys what we call a normal ordering property  and the representation
satisfies a so-called `Popescu condition'. This implies in turn that any such representation admits an isometric dilation. These
sufficient conditions for the existence of an isometric dilation should be compared with recent results of \cite{Solk}, where
sufficient and necessary conditions for the existence of a \emph{regular} dilation were established. The dilations constructed
via the Poisson transform in the case of product systems related to graphs are of a \emph{$^*$-regular} type. It is shown that in
general a minimal (not necessarily regular) isometric dilation of a contractive representation is doubly commuting if and only if
it satisfies the $^*$-regularity property. Contrary to the classical case of commuting Hilbert space contractions, here the
difference between the regular and $^*$-regular dilations is fundamental, as there is no natural adjoint operation on a class of
representations of a given product system (moreover we cannot always assume that the dilations have natural `unitary' extensions,
see \cite{Wold} and references therein).

In the second part of the paper we consider the case of certain families of contractions associated with a higher-rank graph
$\La$ and formalise heuristic observations listed in the second paragraph of this introduction. It is shown that the dilations of
\cite{graphdil} can indeed be viewed as dilations of representations of the canonical product system $\be(\La)$.  General results
of the first part of the paper specialised to this context can be interpreted as giving sufficient conditions on existence of
regular or $^*$-regular dilations of certain tuples of contractions satisfying the commutation relations encoded by a higher-rank
graph. In particular one can deduce immediately from \cite{Soltwo} that any $\La$-contraction associated with a rank-2 graph has
a dilation to a Toeplitz-type family.

The detailed plan of the paper is as follows: after listing some general notations we proceed to introduce in Section 1 basic
notions of $C^*$-correspondences, their product systems over $\Nr$ and (covariant completely) contractive  representations of
such objects. In Section 2 we proceed to define isometric dilations of contractive representations and to quote fundamental
results of B.\,Solel (\cite{Soltwo}, \cite{Solk}) on the existence of dilations in the two-dimensional case and on sufficient and
necessary conditions for the existence of regular dilations. A notion of a $^*$-regular isometric dilation is also introduced and
a fact that a minimal isometric dilation is $^*$-regular if and only if it is doubly commuting established. Section 3 is devoted
to the construction of a generalised Poisson transform associated to a representation of a product system with the normal
ordering property satisfying the Popescu condition and to applications of the transform to isometric dilations. In Section 4 we
recall the canonical product system of $C^*$-correspondences associated to a higher-rank graph $\La$ (\cite{toep}) and describe
its representations in terms of the $\La$-families of operators on a Hilbert space. Finally Section 5 presents the general
results of Sections 2 and 3 specified and adapted to the case of dilations of $\La$-families described in Section 4.


Let $\bn_0= \bn \cup \{0\}$. Fix now and for the rest of the paper $r \in \bn$. The canonical `basis' in $\Nr$ will be denoted by
$(e_1, \ldots, e_r)$, with $e:=\sum_{i=1}^r e_i$. The componentwise maximum (respectively, minimum) of $n,m \in \Zr$ is denoted
by $n \vee m$ (respectively, $n \wedge m$) and we write
 $|n|= n_1 + \cdots + n_r$, $n_+=n \vee 0$, $n_- = - (n \wedge 0)$.

\section{Product systems of $C^*$-correspondences and their representations}\label{prodsys}

Let $\alg$ be a $C^*$-algebra. By a $C^*$-correspondence $E$ over $\alg$ is meant a Hilbert $C^*$-module over $\alg$, equipped
with the structure of a left $\alg$-module (via a nonzero $*$-homomorphism $\phi$ mapping $\alg$ into the $C^*$-algebra of
adjointable operators on $E$). $E$ is \emph{essential} as a left $\alg$-module if the closed linear span of $\phi(\alg) E$ is
equal to $E$. Each $C^*$-correspondence is considered with the usual operator space structure (i.e.\ the one coming from viewing
it as a corner in the appropriate linking algebra). The $C^*$-algebra of adjointable operators on $E$ is denoted by
$\mathcal{L}(E)$. Further details can be found in \cite{Lance} or \cite{Morita}; note  that we will often use the concept of
internal tensor products in the category of Hilbert $C^*$-modules equipped with left actions. In particular any representation
$\sigma$ of $\alg$ on a Hilbert space $\Hil$ allows us to consider a new Hilbert space $E \ot_{\sigma} \Hil$ equipped with the
representation of $\alg$ arising from the left action of $\alg$ on $E$.

Fix now and for the rest of the paper $r \in \bn$. As explained in [So$_{1-2}$] a \emph{product system $\be$ of
$C^*$-correspondences over $\bn_0^r$}
 (formally introduced in \cite{fowl} for a general countable semigroup with a neutral element)
 can be thought of as a family of $r$ $C^*$-correspondences
$\{E_1, \ldots,E_r\}$ over the same $C^*$-algebra together with the unitary isomorphisms $t_{i,j}: E_i \ot E_j \to E_j \ot E_i$
($i>j$) satisfying the natural associativity relations:
\[ (\id_{E_l} \ot t_{i,j}) (t_{i,l} \ot \id_{E_j}) (\id_{E_i} \ot t_{j,l}) =
 (t_{j,l} \ot \id_{E_i}) (\id_{E_j} \ot t_{i,l}) (t_{i,j} \ot \id_{E_l})\]
 for all $1\leq i<j<l\leq k$.
This point of view entails identifying for all $n=(n_1, \ldots, n_r) \in \bn_0^r$ the correspondence $\be (n)$ with $E_1^{\ot^{
n_1}} \ot \cdots \ot E_r^{\ot^{ n_r}}$. We write $t_{i,i} = \id_{E_i \ot E_i}$, $t_{i,j} = t_{j,i}^{-1}$ for $i<j$ and also
define unitary isomorphisms $t_{m,n}:\be(m) \ot \be(n) \to \be(n) \ot \be(m)$ ($m,n \in \Nr$) by obvious compositions of tensor
extensions of appropriate $t_{i,j}$'s.

Let $\Fock_{\be}:= \bigoplus_{n \in \Nr} \be(n)$ denote the Fock module of $\be$ (see \cite{fowl} for the details of the
construction). It is a $C^*$-correspondence over $\alg$. For each $n \in \Nr$ and $e \in \be(n)$ define the creation operator
$L_e: \Fock_{\be} \to \Fock_{\be}$ by the formula
\[ L_e (f) = e \ot f, \;\; f \in \Fock_{\be}.\]
The Toeplitz algebra associated with $\be$ is a concrete $C^*$-algebra in
$\mathcal{L}(\Fock_{\be})$ generated by all creation operators as above. It will be denoted by
$\Toep$.

\begin{deft}
A product system $\be$ of $C^*$-correspondences over $\Nr$ is called compactly aligned if given
$n,m \in \Nr$ and two operators $S\in \mathcal{K}(\be(n))$, $T\in \mathcal{K}(\be(m))$ the operator
$S_n^{n \vee m} T_m^{n \vee m} \in \mathcal{K}(\be(n \vee m))$, where $S_n^{n \vee m}:= S \ot
I_{\be(n\vee m- n)}$ and $T_m^{n \vee m}:= T \ot I_{\be(n\vee m- m)}$.
\end{deft}

The notion of  compact alignment may seem rather technical, but it proved to be very useful (\cite{fowl}). For the product system
associated with a higher-rank graph it is equivalent to the graph in question being \emph{finitely aligned} (see Sections
\ref{graphrep} and \ref{dilgraph}). Examples coming from graphs suggest also that compact alignment of a product system is
closely related to a form of `normal ordering' in the Toeplitz algebra. As we have not been able to determine whether these two
properties coincide in general, we introduce the following definition:

\begin{deft}
A product system $\be$ of $C^*$-correspondences over $\Nr$ is said to have a normal ordering property if
$\Toep=\ol{\textrm{Lin}}\{L_e L_f^*:e,f \in \bigcup_{n \in \Nr}\be(n)\}$.
\end{deft}

The normal ordering property may be thought of as a strong form of `double commutativity' of the creation operators in the
Toeplitz algebra. This is naturally seen when we work with product systems associated with higher-rank graphs in Sections
\ref{graphrep} and \ref{dilgraph}. Note that if $\alg=\bc$ then each $E_j$ is a Hilbert space and the structure of a product
system is determined by the Hilbert space unitaries $t_{i,j}:E_i \to E_j$ (precise description can be found in \cite{Solk} or in
\cite{Wold}). If each $E_j$ is additionally assumed to be finite-dimensional we are in the situation analysed in
\cite{SteveBaruch} and it is easy to see that the corresponding product system has a normal ordering property (and is compactly
aligned).

\subsection*{Representations of $C^*$-correspondences}
The notions presented  in this subsection have been introduced and developed in the series of papers by P.\,Muhly and B.\,Solel
(see \cite{MSgen} and references therein).

\begin{deft}
Let $\Hil$ be a Hilbert space. By a (completely contractive covariant) representation of a
$C^*$-correspondence $E$ over a $C^*$-algebra $\alg$ on $\Hil$ is meant a pair $(\sigma, T)$, where
$(\sigma,\Hil)$ is a representation of $\alg$ on $\Hil$, and $T:E \to B(\Hil)$ is a linear
completely contractive map such that
\[ T(a \xi b) = \sigma(a) T(\xi) \sigma(b), \;\;\; a,b \in \alg, \xi \in E.\]
It is called isometric if for each $\xi, \eta \in E$
\[ T(\xi)^* T(\eta) = \sigma(\langle \xi, \eta \ra).\]
\end{deft}

A representation $(T, \sigma)$ determines a contraction $\wT : E \ot_{\sigma} \Hil \to \Hil$ given
by $\wT (\xi \ot h) = T(\xi) h$ ($\xi \in E, h \in \Hil$). This satisfies:
\begin{equation} \label{tilde} \wT (\phi(a) \ot I_{\Hil}) = \sigma(a) \wT, \;\;\; a \in \alg
\end{equation}
($\phi$ denoting the left action of $\alg$ on $E$), and one can in fact show that, given a
representation $\sigma$, there is a 1-1 correspondence between contractions satisfying
\eqref{tilde} and representations of $E$ (\cite{MSgen} Lemma 2.1). The isometric representations
are exactly those for which $\wT$ is an isometry. The representation $(\sigma, T)$ is said to be
\emph{(fully) coisometric} if $\wt{T}\wt{T}^* = I_{\Hil}$.

It is easy to see how the notion of a representation of a $C^*$-correspondence extends to a product
system.

\begin{deft} \label{prodrep}
Let $\be$ be a product system of $C^*$-correspondences over $\Nr$. By a (completely contractive covariant) representation of
$\be$ on a Hilbert space $\Hil$ is meant a tuple $(\sigma, T^{(1)}, \ldots, T^{(r)})$, where $(\sigma,T^{(i)})$ is a
representation of $E_i$ on $\Hil$ and
\begin{equation} \label{rep} \wT^{(i)} (I_{E_i} \ot \wT^{(j)}) = \wT^{(j)} (I_{E_j} \ot \wT^{(i)}) (t_{i,j} \ot I_{\Hil})\end{equation}
for $i,j\in \oner$. Such a representation is called isometric if each $(\sigma, T^{(i)})$ is isometric, and coisometric if each
$(\sigma, T^{(i)})$ is coisometric.
\end{deft}

To lighten the notation we will occasionally write $\repT$ for $(\sigma, T^{(1)}, \ldots, T^{(r)})$. We will also exploit the
inductively defined maps $T(n)(e)\in B(\Hil)$ ($n \in \Nr$, $e \in \be(n)$) (see  \cite{MSgen}  or \cite{Wold}) and their natural
partners $\wt{T}(n):\be(n) \ot_{\sigma} \Hil \to \Hil$. It is important to note that because of the condition \eqref{tilde}
operators $\wt{T}(n) \wt{T}(n)^*$ belong to $\sigma(\alg)'$.

If we represent the Toeplitz algebra faithfully on a Hilbert space, then the map $\bigcup_{n \in \Nr}\be(n)\ni e \to L_e \in
\Toep$ yields in a natural way a representation of $\be$, called further the Fock-Toeplitz representation. It is easily seen to
be isometric.

In what follows we will often consider doubly commuting representations; these have especially good properties  in terms of the
dilations or Wold decompositions (see respectively \cite{Solk} and \cite{Wold}).

\begin{deft}  \label{dcom}
A representation $(\sigma, T^{(1)}, \ldots, T^{(r)})$ of  $\be$ on a Hilbert space $\Hil$ is called doubly commuting if for each
$i,j \in \{1,\ldots,r\}$, $i \neq j$ implies
\begin{equation}\wT^{(j)^*} \wT^{(i)} =
 (I_{E_j} \ot \wT^{(i)})  (t_{i,j} \ot I_{\Hil})  (I_{E_i} \ot \wT^{(j)^*}).
 \label{doubly}\end{equation}
\end{deft}

For isometric representations of product systems over $\Nr$  double commutativity is exactly the same as  Nica-covariance
considered in \cite{toep} (\cite{Solk}, Remark 3.12). It also has the following equivalent characterisation:

\begin{lem} \label{kerlem}
An isometric representation $(\sigma, T^{(1)}, \ldots, T^{(r)})$ of  $\be$ on a Hilbert space $\Hil$ is doubly commuting if and
only if for each $i,j \in \{1,\ldots,r\}$, $i \neq j$
\begin{equation} \label{kernels} \wt{T}^{(i)} (\Ker (I_{E_i} \ot \wt{T}^{(j)^*})) \subset
\Ker (\wt{T}^{(j)^*}). \end{equation}
\end{lem}

\begin{proof}
Let $i,j$ be as above and denote the operator $I_{E_i} \ot \wt{T}^{(j)}:E_i \ot E_j \ot_{\sigma} \Hil \to E_i \ot_{\sigma} \Hil$
by $\Gamma_{ij}$. Note that the $\Ker (\Gamma_{ij}^*) = E_i \ot_{\sigma} \Ker (\wt{T}^{(j)^*})$. It can be proved exactly in the
same way as the well known statement for kernel of the operator $I_{\Kil_1} \ot S$, where $\Kil_1, \Kil_2$ are Hilbert spaces and
$S \in B(\Kil_2)$ (at least if you know how to show the latter without using an orthonormal basis in $\Kil_1$).

If $\repT$ is doubly commuting, then \eqref{kernels} is easily seen to be satisfied. Suppose then that  \eqref{kernels} holds.
Any vector in $E_i \ot_{\sigma} \Hil$ can be decomposed as a sum of an element in $\Ker (\Gamma_{ij}^*)$ and in $\Ker
(\Gamma_{ij}^*)^{\perp} = \ol{ \Ran(\Gamma_{ij})}$.  It is therefore enough to show that the both sides of the equation
\eqref{doubly} hold on  $\Ran(\Gamma_{ij})$. Let then $z \in E_i \ot E_j \ot_{\sigma} \Hil$. Then
\begin{align*} \wT^{(j)^*} \wT^{(i)}
\Gamma_{ij} (z) &= \wT^{(j)^*} \wT^{(i)} (I_{E_i} \ot \wt{T}^{(j)}) (z) = \wT^{(j)^*} \wT^{(j)} (I_{E_j} \ot \wT^{(i)}) (t_{i,j}
 \ot I_{\Hil})(z) \\&= (I_{E_j} \ot \wT^{(i)}) (t_{i,j} \ot I_{\Hil})(I_{E_i} \ot \wT^{(j)^*} \wT^{(j)} )(z)
\\&= (I_{E_j} \ot \wT^{(i)}) (t_{i,j} \ot I_{\Hil})(I_{E_i} \ot \wT^{(j)^*}) (\Gamma_{ij}(z)).
\end{align*}
This ends the proof.
\end{proof}

\section{General properties of isometric dilations of representations of product systems of $C^*$-correspondences}
\label{dilgen}

In this section we discuss several classes of isometric dilations of a representation of a product system of
 $C^*$-correspondences.

\begin{deft}
Let $\be$ be a product system of $C^*$-correspondences over $\bn_0^r$ and let $(\sigma,
T^{(1)},\ldots, T^{(r)})$ be a representation of $\be$ on a Hilbert space $\Kil$. We say that
$(\pi, V^{(1)},\ldots, V^{(r)})$, an isometric representation of $\be$ on a Hilbert space $\Kil
\supset \Hil$, is an isometric dilation of $(\sigma, T^{(1)},\ldots, T^{(r)})$ if
\begin{rlist}
\item $\forall_{a \in \alg}\;\;\; \pi(a)|_{\Hil} = \sigma(a)$;
\item $\forall_{i \in \oner}\; \forall_{\xi \in \Hil} \;\;   (\wt{V}^{(i)})^* (\xi) = (\wt{T}^{(i)})^* (\xi)$.
\end{rlist}
The dilation is called minimal if $\Kil= \ol{\Lin}\{V(n)(e) \xi: n \in \Nr, e \in \be(n), \xi \in
\Hil\}$.
\end{deft}

Note that the condition (ii) above in particular exploits the identification of $E_i \ot_{\sigma} \Hil$ with  a subspace of $E_i
\ot_{\pi} \Kil$. Moreover  $P_{\Hil} \in \pi(\alg)'$, so that the operators of the form $T \ot P_{\Hil}$ ($T \in
\mathcal{L}(\be(n))$) are well defined as operators in $\mathcal{L}(\be(n) \ot_{\pi} \Kil)$. This is used in Definition
\ref{regular} below.

It is known that two commuting contractions can be always jointly dilated to commuting isometries. The following result
established by B.\,Solel shows that this phenomenon persists in the category of representations of $C^*$-correspondences.

\begin{tw}[[So$_1${]}, Theorem 4.4] \label{sol1}
Let $\be$ be a product system of $C^*$-correspondences over $\bn_0^2$. Every (completely contractive covariant) representation
$(\sigma, T^{(1)}, T^{(2)})$ of $\be$ on a Hilbert space has a minimal isometric dilation $(\pi, V^{(1)}, V^{(2)})$. If $\sigma$
is nondegenerate and $E_1, E_2$ are essential then $\rho$ is nondegenerate.
\end{tw}

\subsection*{Regular isometric dilations (after B.\,Solel)}

To formulate the next result we need a few more definitions. For $u=\{u_1,\ldots, u_k\}  \subset
\oner$ write $e(u) = e_{u_1} + \cdots + e_{u_k}$.

\begin{deft} \label{Solelcond}
Let $\be$ be a product system of $C^*$-correspondences over $\bn_0^r$. A representation $(\sigma, T^{(1)}, \ldots, T^{(r)})$ of
$\be$ on a Hilbert space is said to satisfy the Brehmer-Solel condition if for each $v \subset \oner$
\[ \sum_{u\subset v} (-1)^{|u|} (I_{\be(e(v) - e(u))} \ot \wt{T}(e(u))^* \wt{T}(e(u))) \geq 0.\]
\end{deft}

The condition above first appeared in the context of commuting families of contractions in \cite{Brehmer}; recently it was
exploited in the context of product systems of $C^*$-correspondences in \cite{Solk}.

\begin{deft} \label{regular}
Let $(\sigma, T^{(1)}, \ldots, T^{(r)})$ be a representation of $\be$ on a Hilbert space $\Hil$. An isometric dilation $\repV$ of
$\repT$ is said to be regular if for all $n \in \Zr$
\begin{equation}\label{reg} (I_{\be(n_{-})} \ot P_{\Hil}) \wt{V}(n_{-})^*
\wt{V}(n_+)|_{\be(n_+) \ot_{\sigma} \Hil} = \wt{T}(n_{-})^* \wt{T}(n_+).\end{equation}
\end{deft}

B.\,Solel showed that the condition described in Definition \ref{Solelcond} characterises these representations which allow
\emph{regular} isometric dilations.

\begin{tw}[[So$_2${]}, Theorem 3.5] \label{sol2}
Let  $\be$ be a product system of $C^*$-correspondences over $\bn_0^r$. A (completely contractive covariant) representation
$(\sigma, T^{(1)},\ldots, T^{(r)})$ of $\be$ on a Hilbert space has a regular isometric dilation if and only if it satisfies the
Brehmer-Solel condition.
\end{tw}

Note that minimal regular dilations are necessarily unique, in the sense that any two such
dilations respectively on Hilbert spaces $\Kil$ and $\Kil'$ are intertwined by a unitary $U: \Kil
\to \Kil'$.

\subsection*{$^*$-regular dilations}

Let us begin with a simple equivalent characterisation of regularity of an isometric dilation.

\begin{lem} \label{simpreg}
An isometric dilation $(\pi, V^{(1)},\ldots, V^{(r)})$ of a representation \\$(\sigma,
T^{(1)},\ldots, T^{(r)})$  of $\be$ is regular if and only if for all $n,m \in \Nr$ such that $n_j
\neq 0$ implies $m_j= 0$ ($j \in \oner$) and all $e\in \be(n)$, $f \in \be(m)$
\begin{equation} \label{regsimp} P_{\Hil} (V(n)(e))^* V(m)(f)|_{\Hil} = (T(n)(e))^* T(m)(f).\end{equation}
\end{lem}

\begin{proof}
Note that condition \eqref{reg} is satisfied for all $n \in \oner$ if and only if for all $n,m \in
\Nr$ such that $n_j \neq 0$ implies $m_j= 0$ ($j \in \oner$) there is
\[(I_{\be(n)} \ot P_{\Hil}) \wt{V}(n)^* \wt{V}(m)|_{\be(m) \ot \Hil} =
\wt{T}(n)^* \wt{T}(m).\] The last condition is equivalent to the fact that for all $\xi\in \Hil$, $f \in \be(m)$
\[ (I_{\be(n)} \ot P_{\Hil}) \wt{V}(n)^* V(m)(f) \xi =
\wt{T}(n)^* T(m) (f) \xi,\] and further to the fact that for all $e \in \be(n), \eta \in \Hil$
\[ \langle e \ot \eta, (I_{\be(n)} \ot P_{\Hil}) \wt{V}(n)^* V(m)(f) \xi \rangle = \langle e \ot \eta, \wt{T}(n)^* T(m) (f)
\xi \rangle.\] This in turn holds if and only if
\[ \langle V(n) (e) \eta, V(m)(f) \xi  \rangle  = \langle T(n) (e) \eta, T(m)(f) \xi  \rangle,\]
if and only if
\[ \langle  \eta, P_{\Hil} (V(n) (e))^*V(m)(f) \xi  \rangle  = \langle  \eta, (T(n) (e))^*T(m)(f) \xi  \rangle.\]
This ends the proof.
\end{proof}

Recall (\cite{Timot}) that if $(S,T)$ is a commuting pair of contractions on a Hilbert space $\Hil$ and $(U,V)$ is a commuting
isometric dilation of $(S,T)$ then it is said to be $^*$-regular if for all $k,l \in \bn$
\[P_{\Hil} (U^*)^l V^k|_{\Hil} = T^k (S^*)^l.  \]

How should a corresponding definition look here? Note that if $\repV=(\pi, V^{(1)},\ldots, V^{(r)})$ is an isometric dilation  of
a representation $\repT=(\sigma, T^{(1)},\ldots, T^{(r)})$  of $\be$ then for all $n,m \in \Nr$ and all $e\in \be(n)$, $ f \in
\be(m)$
\begin{equation} \label{correct0} P_{\Hil} V(n)(e) (V(m)(f))^*|_{\Hil} = T(n)(e) (T(m)(f))^*.\end{equation}
 Moreover one can see that similarly for all $n \in
\Zr$
\begin{align} \label{correct} (I_{\be(n_{-})} \ot P_{\Hil})  &(I_{\be(n_{-})} \ot \wt{V}(n_+)) (t_{n_{+},
n_{-}} \ot I_{\Kil}) (I_{\be(n_+)} \ot \wt{V}(n_{-}))^*)|_{\be(n_+) \ot \Hil} \\ \notag =&(I_{\be(n_{-})} \ot \wt{T}(n_+))
(t_{n_{+}, n_{-}} \ot I_{\Hil}) (I_{\be(n_+)} \ot \wt{T}(n_{-}))^*);\end{align} it is enough to observe that from the definition
of the dilation it follows that for all $n,m \in \Nr$
\[ (I_{\be(n)} \ot \wt{V}(m))^*|_{\be(n) \ot_{\sigma} \Hil} = (I_{\be(n)} \ot \wt{T}(m))^*|_{\be(n) \ot_{\sigma} \Hil}.\]

If the dilation $\repV$ is doubly commuting then the condition \eqref{correct} reduces to \[(I_{\be(n_{-})} \ot P_{\Hil})
\wt{V}(n_{-})^* \wt{V}(n_+) |_{\be(n_+) \ot_{\sigma} \Hil} = (I_{\be(n_{-})} \ot \wt{T}(n_+)) (t_{n_{+}, n_{-}} \ot I_{\Hil})
(I_{\be(n_+)} \ot \wt{T}(n_{-}))^*.\] The latter can be seen as a natural generalisation of the notion of $^*$-regularity.

In the classical context of commuting contractions Theorem 2 of \cite{Timot} (see also \cite{GaSu})
shows that a minimal isometric dilation is $^*$-regular if and only if it is doubly commuting. The
same remains true in our context, as the next theorem shows. The proof is a natural generalisation
of that in \cite{Timot}. The basic idea is the following: as the last equation implies  `double
commutativity on $\Hil$', we need to exploit minimality to deduce `double commutativity on $\Kil$'.

\begin{tw} \label{dc=streg}
A minimal isometric dilation $(\pi, V^{(1)},\ldots, V^{(r)})$ of a representation $(\sigma, T^{(1)},\ldots, T^{(r)})$  of $\be$
is doubly commuting if and only if it is $^*$-regular, that is if for all $n \in \Zr$
 \begin{align} \label{streg}(I_{\be(n_{-})} \ot P_{\Hil})&  \wt{V}(n_{-})^* \wt{V}(n_+)  |_{\be(n_+) \ot_{\sigma} \Hil}
= \\& \notag (I_{\be(n_{-})} \ot \wt{T}(n_+)) (t_{n_{+}, n_{-}} \ot I_{\Hil}) (I_{\be(n_+)} \ot \wt{T}(n_{-}))^*.\end{align}
\end{tw}

\begin{proof}
The fact that a doubly commuting dilation is automatically $^*$-regular has been explained in the discussion before the theorem
(and does not require minimality). Suppose then that $\repV$ (acting on $\Kil$) is a minimal $^*$-regular dilation of $\repT$
(acting on $\Hil$) and fix $i,j \in \oner$, $i \neq j$. By Lemma \ref{kerlem} it is enough to show that $\wt{V}^{(i)} (\Ker
(\Gamma_{ij}^*)) \subset \Ker (\wt{V}^{(j)^*})$, where $\Gamma_{ij}= I_{E_i} \ot \wt{V}^{(j)}$. As $\repV$ is minimal, $E_i \ot
E_j \ot_{\pi} \Kil$ is generated by $\{e \ot f  \ot \wt{V}(m) (g \ot  \xi): e \in E_i, f \in E_j, m \in \Nr, g \in \be(m), \xi
\in \Hil\}$. This implies that $\Ran (\Gamma_{ij})$ is generated by
\[\{e \ot \wt{V}^{(j)}( f \ot\wt{V}(m) (g \ot \xi)) : e \in
E_i, f \in  E_j, m \in \Nr, g \in \be(m),\xi \in \Hil\},\] so also by
\[\{e \ot  \wt{V}(m) (g \ot \xi)  : e \in
E_i,  m \in \Nr, m_j \neq 0, g \in \be(m),\xi \in \Hil\}.\] This  in turn implies that $\Kil_0:=\Ker(\Gamma_{ij}^*) =
\Ran(\Gamma_{ij})^{\perp} \subset E_i \ot_{\pi} \Kil$ is equal to the subspace generated by
\[\left\{P_{\Kil_0} \left(e \ot  \wt{V}(m) (g \ot \xi)\right)  : e \in
E_i,  m \in \Nr, m_j = 0, g \in \be(m),\xi \in \Hil\right\}.\]
 Suppose then that $e \in E_i,  m \in \Nr, m_j = 0, g \in \be(m),\xi \in \Hil$.
We will show that
\begin{equation} P_{\Kil_0} (e \ot  \wt{V}(m) (g \ot \xi) ) = e \ot \wt{V}(m) (g \ot \xi - g \ot \wt{V}^{(j)}
\wt{T}^{(j)^*} \xi). \label{projection} \end{equation}
Note that as $e$ does not play any significant role here, remembering the
remarks made in the beginning of the proof of Lemma \ref{kerlem} it is enough to show that
\[ P_{\Kil_1} z = z - v,\]
where $ z=\wt{V}(m) (g \ot \xi)$, $ v = \wt{V}(m) (g \ot \wt{V}^{(j)} \wt{T}^{(j)^*} \xi)$ and $\Kil_1 = \Ker (\wt{V}^{(j)^*})$.
As $v \in \Ran(\wt{V}^{(j)})$, it suffices if we can show that $z-v \perp \Ran(\wt{V}^{(j)})$. We are going to exploit minimality
once more. Let $f \in E_j, n \in \Nr$, $h \in \be(n), \eta \in \Hil$ and compute
\[
A :=\langle \wt{V}^{(j)} (f \ot \wt{V}(n) (h \ot \eta)), z \rangle = \langle  f \ot h \ot \eta, \wt{V}(n+e_j)^* \wt{V} (m) (g \ot
\xi)\rangle\] Let now $l = (n+e_j) \wedge m$, $p = n+e_j -l$, $q=m-l$. Note that $l_j = 0$, $p_j \neq 0$. Then
\begin{align*} A &= \langle  f \ot
h \ot \eta, (t_{l,j} \ot I_{\be(p-e_j)} \ot P_{\Hil}) ( I_{\be(l)} \ot \wt{V} (p)^*) \wt{V}(l)^*\wt{V}(l) ( I_{\be(l)} \ot \wt{V}
(q)) (g \ot \xi)\rangle \\ &= \langle  (t_{j,l} \ot I_{\be(p-e_j)\ot_{\sigma} \Hil}) (f \ot h \ot \eta),  (I_{\be(l)} \ot
I_{\be(p)} \ot P_{\Hil}) ( I_{\be(l)} \ot \wt{V} (p)^*) ( I_{\be(l)} \ot \wt{V} (q))  (g \ot \xi)\rangle\end{align*} The
$^*$-regularity condition \eqref{streg} implies that for all $l \in \Nr$ and all $p,q \in \Nr$ such that $p \wedge q =0$ there is
\begin{align*} (I_{\be(l)} & \ot  I_{\be(p)} \ot P_{\Hil}) (I_{\be(l)} \ot  \wt{V}(p)^*) (I_{\be(l)} \ot \wt{V}(q)) |_{\be(l+q) \ot_{\sigma} \Hil}
\\&= (I_{\be(l)} \ot I_{\be(p)} \ot \wt{T}(q)) (I_{\be(l)} \ot t_{q,p} \ot I_{\Hil}) (I_{\be(l)} \ot I_{\be(q)} \ot
\wt{T}(p))^*)
\end{align*}
Therefore
\begin{align*} A =&  \langle  (t_{j,l} \ot I_{\be(p-e_j)\ot_{\sigma} \Hil}) (f \ot h \ot \eta),
\\&
 (I_{\be(l)} \ot I_{\be(p)} \ot \wt{T}(q)) (I_{\be(l)} \ot t_{q,p} \ot I_{\Hil}) (I_{\be(l)} \ot I_{\be(q)} \ot \wt{T}(p)^*) (g
\ot \xi)\rangle.
\end{align*}
Similarly
\begin{align*}
B:= & \langle \wt{V}^{(j)} (f \ot \wt{V}(n) (h \ot \eta)), v \rangle = \langle  f \ot h \ot \eta, \wt{V}(n+e_j)^* \wt{V} (m) (g
\ot
 \wt{V}^{(j)} \wt{T}^{(j)^*} \xi )\rangle\\& =
\langle  f \ot h \ot \eta, \wt{V}(n+e_j)^* \wt{V} (m+e_j) (g \ot \wt{T}^{(j)^*} \xi )\rangle.
\end{align*}
Put now $l' =  (n+e_j) \wedge (m+e_j)$, $p' = n+e_j -l'$, $q'=m+e_j-l'$. Note that $l'=l+e_j$, $p'=p-e_j$, $ q'= q$. Continuing
as before we obtain (note that $t_{j,l}$ no longer features, as $l'_j \neq 0$)
\begin{align*} B &= \langle  f \ot
h \ot \eta, (I_{l'} \ot I_{\be(p')} \ot P_{\Hil}) ( I_{\be(l')} \ot \wt{V} (p')^*) \wt{V}(l')^*\wt{V}(l') ( I_{\be(l')} \ot
\wt{V} (q')) (g \ot \wt{T}^{(j)^*} \xi)\rangle
\\ &= \langle  f \ot h \ot \eta,
 (I_{\be(l')} \ot I_{\be(p')} \ot \wt{T}(q')) (I_{\be(l')} \ot t_{q',p'} \ot I_{\Hil}) (I_{\be(l')} \ot I_{\be(q')} \ot \wt{T}(p')^*) (g
\ot \wt{T}^{(j)^*} \xi)\rangle
\\&= \langle  f \ot h \ot \eta,  (I_{\be(l+e_j)} \ot I_{\be(p-e_j)} \ot \wt{T}(q)) \\
&  \;\;\hspace*{2cm} (I_{\be(l+e_j)} \ot t_{q,p-e_j} \ot I_{\Hil}) (I_{\be(l)} \ot t_{q,j} \ot I_{\be(p-e_j)})
 (I_{\be(l)} \ot I_{\be(q)} \ot \wt{T}(p)^*) (g
\ot  \xi)\rangle.
\end{align*}
The comparison of the formulas above shows that $A=B$ if only
\begin{align*} (t_{l,j} & \ot I_{\be(p-e_j)}\ot I_{\be(q)} ) (I_{\be(l)} \ot
I_{\be(p)} \ot I_{\be(q)}) (I_{\be(l)} \ot t_{q,p} ) (I_{\be(l)} \ot I_{\be(q)} \ot I_{\be(p)})
\\&= (I_{\be(l+e_j)} \ot
I_{\be(p-e_j)} \ot I_{\be(q)}) (I_{\be(l+e_j)} \ot t_{q,p-e_j}) (I_{\be(l)} \ot t_{q,j} \ot I_{\be(p-e_j)}) (I_{\be(l)} \ot
I_{\be(q)} \ot I_{\be(p)})
\end{align*}
This can be further reduced to checking two equalities
\[t_{l,j}  \ot I_{\be(p-e_j)} =I_{\be(l+e_j)} \ot
I_{\be(p-e_j)} \] and
\[ I_{\be(l)} \ot t_{q,p}  = (I_{\be(l+e_j)} \ot t_{q,p-e_j}) (I_{\be(l)} \ot t_{q,j} \ot I_{\be(p-e_j)});\]
these finally are simple consequences of the definition of $t_{m,n}$ in the beginning of Section \ref{prodsys}.

The equality $A=B$ implies that
\[ \langle \wt{V}^{(j)} (f \ot \wt{V}(n) (h \ot \eta)), z-v \rangle = 0.\]
As $f \in E_j, n \in \Nr$, $h \in \be(n), \eta \in \Hil$ are arbitrary and $\repV$ is minimal, $z-v \perp \Ran(\wt{V}^{(j)})$ and
\eqref{projection} is proved. Note now that
\begin{align*} \wt{V}^{(i)} ( P_{\Kil_0} (e \ot z)) &=  \wt{V}^{(i)}  (e\ot (z - v))\\&=
\wt{V}^{(i)} (e \ot \wt{V}(m) (g \ot \xi - g \ot \wt{V}^{(j)} \wt{T}^{(j)^*} \xi)) \\&=  \wt{V}^{(m+e_i)} (e \ot g \ot \xi - e
\ot g \ot \wt{V}^{(j)} \wt{T}^{(j)^*} \xi) \\&= P_{\Kil_0} \wt{V}^{(m+e_i)} (e \ot g \ot \xi) \in \Kil_0.\end{align*}
 This ends the
proof.
\end{proof}

It follows from the theorem above that a minimal isometric doubly commuting dilation of a
representation of a product system is unique up to a unitary equivalence, as condition
\eqref{streg} together with minimality determines scalar products between all vectors in $\Kil$. In
general a minimal isometric dilation need not be unique. Concrete examples of this phenomenon can
be found in \cite{DPY} (Examples 4.3 and 4.4).

The following result can be shown in a similar way to Proposition 2.6 in \cite{Wold}.

\begin{prop}
If $\repV$ is a minimal isometric doubly commuting dilation of a coisometric representation
$\repT$, then $\repV$ is coisometric.
\end{prop}

\section{Generalised Poisson transform and isometric dilations}

In this section we describe how to construct isometric dilations via the generalised Poisson
transform associated with a given representation of a product system. In a similar
multi-dimensional context it has been first introduced in \cite{PPois}; the one-dimensional
counterpart for representations of $W^*$-correspondences has been recently investigated  in
\cite{MSker}.

The next definition describes a natural variation on the type of conditions considered when one wants to construct isometric
dilations of higher-rank objects (see \cite{Nagy} and references therein) and should be compared to the introduced earlier
Brehmer-Solel condition.

\begin{deft} \label{PoTcond}
Let $\be$ be a product system of $C^*$-correspondences over $\bn_0^r$. For a representation
$(\sigma, T^{(1)}, \ldots, T^{(r)})$ of $\be$ on a Hilbert space $\Hil$ define the defect operator
($s \in (0,1)$)
\begin{equation}\label{DeltaT} \PoT =
\sum_{n \in \Nr, n\leq e} (-s^2)^{|n|} \wt{T}(n) \wt{T}(n)^*.
\end{equation}
The representation $\repT$ is said to satisfy the Popescu condition (or condition `P') if there
exists $\rho\in (0,1)$ such that for all $s \in (\rho,1)$ the operator $\PoT$ is positive.
\end{deft}

The condition above in a similar form first appeared in \cite{PPois}. Its variant for families of
contractions associated with higher-rank graphs was extensively studied in \cite{graphdil}. Because
of the condition \eqref{tilde} the defect operator $\PoT$ is in $\sigma(\alg)'$. It is easy to see
that if $\repT$ is doubly commuting or coisometric then it satisfies the Popescu condition.

\begin{tw}\label{PoisE}
Let $\be$ be a product system of $C^*$-correspondences over $\bn_0^r$ having a normal ordering
property and let $\repT=(\sigma, T^{(1)}, \ldots, T^{(r)})$ be a representation of $\be$ on a
Hilbert space $\Hil$ satisfying the Popescu condition. Then there exists a unique continuous linear
map $R_{\repT}:\Toep \to B(\Hil)$ satisfying
\[R_{\repT}(L_{e} L_{f}^*) = T(n)(e) (T(m)(f))^*, \;\;\;\; n,m \in \Nr, e\in \be(n),f \in \be(m).\]
The map $R_{\repT}$ will be called the generalised $\be$-Poisson transform (associated with
$\repT$). It is completely positive and contractive, unital if $\Toep$ is unital.
\end{tw}

\begin{proof}
The proof is almost identical to the one given for the case of $\La$-Poisson transforms associated
with higher-rank graphs in \cite{graphdil}. We will therefore only indicate the main points and
extra difficulties arising here. Let $s \in (0,1)$ and consider the operator $\GaT\in B(\Hil)$
given by
\[ \GaT (\xi) = \sum_{n \in \Nr} s^{2|n|} \wt{T}(n) (I_{\be(n)} \ot \PoT) \wt{T}(n)^*  \xi\]
($\xi \in \Hil$). It can be checked that $\GaT=I_{\Hil}$ (see Lemma 2.1 of \cite{graphdil}). CHECK!!! Let $\rho\in (0,1)$ be such
that for all $s \in (\rho,1)$ the operator $\PoT$ is positive. As $\PoT \in \sigma(\alg)'$, also $(\PoT)^{\frac{1}{2}} \in
\sigma(\alg)'$ and moreover for each $n \in \Nr$ the operator $I_{\be(n)} \ot \PoT$ on $\be(n) \ot_{\sigma} \Hil$ is positive,
\[(I_{\be(n)} \ot
\PoT)^{\frac{1}{2}} = I_{\be(n)} \ot \PoT^{\frac{1}{2}}.\] Similarly, if $n \in \bn^r, T \in \mathcal{L}(\be(n))$, $S \in
\sigma(\alg)'$, then
\[ (T \ot S)^* = T^* \ot S^*.\]
These properties will be further used without any comments.

Define the isometry $\WoT: \Hil \to \Fock_{\be} \ot_{\sigma} \Hil$ by
\[ \WoT \xi =  \bigoplus_{n \in \Nr} s^{|n|} (I_{\be(n)} \ot \PoT^{\frac{1}{2}}) \wt{T}(n)^* \xi.\]
Let the  map $\sTPois:\mathcal{L}(\Fock_{\be}) \to B(\Hil)$ be given by the formula
\[\sTPois (x) = \WoT^* (x \ot I_{\Hil}) \WoT,\]
 It is clear that $\sTPois$ is completely positive and contractive. Moreover for any $e \in \be(n), f \in \be(m)$ ($n,m \in \Nr$)
\begin{equation} \label{TPois}
\sTPois (L_e L_f^*) = s^{|n| + |m|}T(n)(e) T(m)(f)^*.
\end{equation}
Indeed, let $e,f$ be as above. Note first that for all $n'\in \Nr$
\[\wt{T}(n') (L_e \ot I_{\be(n'-n)\ot_{\sigma} \Hil}) =  T(n) (e) \wt{T}(n'-n), \]
so also
\[ (L_e^* \ot I_{\be(n'-n)\ot_{\sigma}\Hil}) \wt{T}(n')^* = \wt{T}(n'-n)^* (T(n) (e))^*.\]
Compute further ($ \xi, \eta \in \Hil$):
\begin{align*}
\lla \eta, \right.& \left. \sTPois (L_{e} L_{f}^*) \xi \rra = \lla \WoT
\eta,  \left( L_{e} L_{f}^* \ot I_{\Hil}  \right)\WoT \xi \rra \\
&= \lla \sum_{n' \in \Nr} s^{|n'|} (I_{\be(n')} \ot \PoT^{\frac{1}{2}}) \wt{T}(n')^*  \eta, \right.\\&\;\;\;\;\;\;\; \left.
\left(L_{e} L_{f}^* \ot I_{\Hil}
\right) \sum_{m' \in \Nr} s^{|m'|} (I_{\be(m')} \ot \PoT^{\frac{1}{2}}) \wt{T}(m')^* \xi \rra \\
=& \lla \sum_{n' \in \Nr, n' \geq n} s^{|n'|} (I_{\be(n'-n)} \ot \PoT^{\frac{1}{2}}) (L_{e}^*  \ot I_{\be(n'-n)\ot_{\sigma}\Hil})
\wt{T}(n')^*  \eta, \right. \\ & \left.  \;\;\;\;\sum_{m' \in \Nr, m' \geq m} s^{|m'|} (I_{\be(m'-m)} \ot \PoT^{\frac{1}{2}})
(L_{f}^* \ot I_{\be(m'-m)\ot_{\sigma}\Hil}) \wt{T}(m')^* \xi \rra \\
=&\sum_{p \in \Nr} \lla s^{|p|+|n|} (I_{\be(p)} \ot \PoT^{\frac{1}{2}}) \wt{T}(p)^* (T(n) (e))^* \xi, \right. \\&
\;\;\;\;\;\;\;\;\; \left. s^{|p|+|m|} (I_{\be(p)} \ot
\PoT^{\frac{1}{2}}) \wt{T}(p)^* (T(m) (f))^* \eta \rra \\
=& s^{|n|+|m|} \lla (T(n) (e))^* \xi,  \left( \sum_{p \in \Nr} s^{2|p|} \wt{T}(p) (I_{\be(p)} \ot \PoT) \wt{T}(p)^* \right) (T(m) (f))^* \eta \rra \\
=& s^{|n|+|m|} \lla (T(n) (e))^* \xi, \GaT (T(m) (f))^* \eta \rra = \lla  \xi, s^{|n| + |m|}T(n)(e) T(m)(f)^* \eta \rra.
\end{align*}
It is now easy to see that by the normal ordering property the limit $\lim_{s \to 1^-} \sTPois (x)$ in the norm topology exists
for each $x \in \Toep$ , and moreover the map $\TPois:\Toep \to B(\Hil)$ defined by
\[ \TPois(x) = \lim_{s
\to 1^-} \sTPois (x), \;\;\; x \in \Toep\] satisfies all the requirements of the theorem. Uniqueness is another consequence of
the normal ordering property.
\end{proof}

\begin{tw} \label{main}
Let $\be$ be a product system of $C^*$-correspondences over $\bn_0^r$ having a normal ordering property and let $\repT=(\sigma,
T^{(1)}, \ldots, T^{(r)})$ be a representation of $\be$ on a Hilbert space $\Hil$ satisfying the Popescu
condition. Then $\repT$ has an 
isometric dilation. 
\end{tw}

\begin{proof}
Consider the minimal Stinespring dilation of the Poisson transform $\TPois$ constructed in Theorem \ref{PoisE}. This provides us
with a Hilbert space $\Kil$, a representation $\rho: \Toep \to B(\Kil)$ and an operator $V\in B(\Hil;\Kil)$ such that for all $x
\in \Toep$
\[ \Pois(x) = V^* \rho(x) V\]
and $\Kil = \overline{\Lin} \{\rho(x)V\xi: x \in \Toep, \xi \in \Hil\}$. We may assume that $V$ is an isometry, if necessary
extending $\TPois$ in the unital manner to the unitisation of $\Toep$ in $\mathcal{L}(\Fock_{\be})$. This allows us to view
$\Hil$ as a subspace of $\Kil$. Define for each $i \in \Nr$, $e \in E_i$
\[ V^{(i)}(e) = \rho(L_e);\]
and for $a \in \alg$ (note that in our framework $\alg= \be(0) \subset \Fock_{\be}$)
\[ \pi(a) = \rho(L_a).\]
It is clear that the tuple $\repV=(\pi, V^{(1)},\ldots, V^{(r)})$ is an isometric representation of $\be$, as it is a
$^*$-homomorphic image of the Fock-Toeplitz representation.

The fact that condition  \eqref{correct0} is satisfied follows directly from the definition of $\repV$, so it remains to
establish that each $(V(n)(e))^*$ ($n \in\Nr$, $e \in \be(n)$) leaves $\Hil$ invariant. By the minimality of the Stinespring
dilation we know that
\[\Kil=\ol{\Lin}\{V(m) (f) (V(p)(g))^* \xi: m,p \in \Nr, f \in \be(m), g \in \be(p), \xi \in \Hil\}.\]
Further given $m,p \in \Nr, f \in \be(m), g \in \be(p)$ and $\xi, \eta \in \Hil$,
\begin{align*} \lla   \right.& \left.  V(m)(f)(V(p)(g))^* \xi,(V(n) (e))^* \eta \rra = \lla V(n) (e) V(m) (f) (V(p)(g))^*\xi, \eta \rra \\&=
\lla V(n+m) (e\ot f)  (V(p)(g))^* \xi, \eta \rra  = \lla P_{\Hil} V(n+m) (e\ot f) (V(p)(g))^* P_{\Hil}\xi, \eta \rra \\&= \lla
T(n+m) (e\ot f) (T(p)(g))^* \xi, \eta \rra  = \lla T(m) (f) (T(p)(g))^* \xi, (T(n) (e))^* \eta \rra \\&= \lla P_{\Hil} T(m) (f)
(T(p)(g))^* P_{\Hil} \xi,  (T(n) (e))^* \eta \rra   \\&= \lla V(m) (f)  (V(p)(g))^*\xi,  (T(n) (e))^* \eta \rra .  \end{align*}
This shows that $(V(n)(e))^*|_{\Hil} = (T(n) (e))^*$. In particular
\[\Kil=\ol{\Lin}\{V(m) (f) \xi: m \in \Nr, f \in \be(m), \xi \in \Hil\}.\]
\end{proof}

The approach via a Poisson transform suggests that the constructed dilation should be $^*$-regular.
If the creation operators in $\Toep$ satisfy some variant of the double commutativity, this will be
the case (see Theorem \ref{quoted}). Recall that minimal $^*$-regular dilations are unique, as
explained in the comments after Theorem \ref{dc=streg}. Once again we see here potential analogies
between the normal ordering property, compact alignment and double commutativity of the creation
operators.

\section{Product system of Hilbert bimodules  associated to a higher rank graph}
\label{graphrep}

In this section we recall a construction of a product system of $C^*$-correspondences associated to
a higher-rank graph $\La$ introduced in \cite{toep} and describe its representations in terms of
the $\La$-families of operators on a Hilbert space.

 A rank-$r$ graph $\La$ is a small category
with set of objects $\Lambda^0$ and shape functor $\sigma : \Lambda \to {\bf N}^r$ (where ${\bf
N}^r$ is viewed as the  category with one object and morphisms $\Nr $) satisfying the
\emph{factorisation property} defined in \cite{kupa}. If $n \in \Nr$ the set of morphisms in $\La$
of shape $n$ is denoted by $\La^n$. Further for each $a \in \La^0$ and $n \in \Nr$ write
$\La^n_a:=\{ \la \in \La: s(\la) = a, \sigma(\la) =n\}$ and  $|\lambda |=|\sigma(\lambda)|$. The
morphisms in $\Lambda$ may be thought of as paths in a `multi-coloured' graph with vertices indexed
by the set $\La^0$. The range and source maps are respectively denoted by $r:\La \to \La^0$ and
$s:\La \to \La^0$. The factorisation property says that if $m,n \in \Nr$ then every morphism
$\lambda \in \Lambda^{m+n}$ is a unique product $\lambda=\mu \nu$ of a  $\mu \in \Lambda^m$ and
$\nu \in \Lambda^n$, where $s(\mu)=r(\nu)$.

A rank-$r$ graph $\La$ is called \emph{finitely aligned} if for each $\la,\mu \in \La$ the set of
minimal common extensions of $\la$ and $\mu$, that is $MCE(\la, \mu):=\{ \nu \in \La :
\exists_{\alpha ,\beta \in \Lambda} \; \nu = \la \alpha = \mu \beta, \: \sigma(\la \alpha) =
\sigma(\la) \vee \sigma(\mu)\}$, is finite.

In \cite{toep}  it was shown that every higher-rank graph can be viewed as a product system of rank-1 graphs and this point of
view leads to associating to such a graph a product system of $C^*$-correspondences. We rephrase this construction below - note
that our conventions on the rank and source follow \cite{book} rather than \cite{toep} and  we are solely interested in product
systems over $\bn_0^r$ (which leads to certain simplifications).

Let $\alg_0= C_0 (\La^0)$ denote the $C^*$-algebra of all complex-valued functions on $\La^0$ vanishing at infinity. Let $j \in
\oner$. Define the $C^*$-correspondence $E_j (\La)$ over $\alg_0$ as follows: $E_j(\La)$ consists of these functions
$x:\La^{e_j}\to\bc$ which are `locally square integrable', i.e.\ for each $a \in \La^0$
\[ x_a:=\sum_{\la \in \La^{e_j}_a} |x(\la)|^2 < \infty \]
and the function $a \to x_a$ vanishes at infinity. The actions of $\alg_0$ on $E_j(\La)$ are
defined via ($f \in \alg_0, x \in E_j(\La), \la \in \La^{e_j}$)
\begin{equation}  (x  \cdot f)
(\la) = x (\la) f(s(\la)), \;\;\;(f \cdot x) (\la) = f(r(\la)) x(\la), \label{action}
\end{equation} and the $\alg_0$ valued scalar product by ($x,y \in E_j(\La), a \in \La^0$)
\begin{equation} \langle x,
y \ra (a) = \sum_{\la \in \La^{e_j}_a} \ol{x(\la)} y (\la).\label{prod}\end{equation} As  finitely supported functions are dense
in $E_j(\La)$, it is easy to see that each of the $C^*$-correspondences $E_j(\La)$ is essential: $\ol{\alg_0 E_j(\La)}=
E_j(\La)$.

It will also be important at a certain point to consider the natural operator space structure of $E_j(\La)$. Intuitively one
should think of $E_j$ as a bundle of Hilbert spaces over $\La^0$ and observe that in the linking algebra picture the Hilbert
spaces in question act as columns. Therefore the natural operator space structure on $E_j$ is the one coming from viewing it as a
bundle of operator spaces $(\Hil_a)_{\textrm{c}}$. The explicit formula for the matricial norms is as follows:
\begin{equation} \label{opstruct} \left\|(x_{ij})_{i,j=1}^n
\right\| = \sup_{a\in \La^0} \left\| \left (\sum_{l=1}^n \sum_{\la \in \La^{e_j}_a}
\ol{x_{li}(\la)} x_{lk}(\la)\right)_{i,k=1}^n \right\|_{M_n}.\end{equation}

 To introduce on $(E_1(\La), \ldots, E_r(\La))$ the structure of a product system
  we identify $E_i(\La) \ot E_j (\La)$ with the space of all functions $z: \La^{e_i + e_j}\to
\bc$ such that  for each $a \in \La^0$
\[ z_a:=\sum_{\la \in \La^{e_i+ e_j}_a} |z(\la)|^2 < \infty \]
and the function $a \to z_a$ vanishes at infinity. The identification is implemented via the
factorisation property: given $\nu \in \La^{e_i+e_j}$ we can decompose it uniquely as $\nu = \la
\mu$, $\nu \in \La^{e_i}, \mu \in \La^{e_j}$ and for $x \in E_i(\La)$, $y \in E_j(\La)$ define
\[ (x \ot y) (\nu):= x(\la) y (\mu).\]
In other words, if $x \in E_i(\La)$, $y \in E_j(\La)$ then for $\nu \in \La^{e_i +e_j}$
\[t_{i,j} (x \ot y) (\nu) = x(\la) y (\mu), \;\; \textrm{where} \;\; \la \in \La^{e_i}, \mu \in \La^{e_j}, \nu = \mu
\la. \] Note that this leads to natural identifications of $\be (n)$ with the spaces of `locally square integrable' functions on
$\La^n$. Precisely speaking, if for each $n \in \Nr$ we define $E_n(\la)$ to be the space of all functions  $x:\La^n \to \bc$
such that for each $a \in \La^0$
\[ x_a:=\sum_{\la \in \La^{n}_a} |x(\la)|^2 < \infty \]
and the function $a \to x_a$ vanishes at infinity. The actions of $\alg_0$ and the $\alg_0$ valued scalar product on $E_n(\La)$
are defined again via formulas \eqref{action} and \eqref{prod} (this time $f \in \alg_0, x,y \in E_n, \la \in \La^n$). Define for
all $n,m \in \oner$ the map $U_{n,m}: E_n(\La) \ot E_n(\La) \to E_{n+m}(\La)$ via the continuous linear extension of the formula
\[ U_{n,m} (x \ot y) (\la)= x(\mu ) y (\nu), \]
where $x \in E_n(\La), y \in E_m(\La), \la \in \La^{n+m}, \la = \mu \nu, \mu \in \La^{n}, \nu \in \La^m$. It can be checked that
$U_{n,m}$ is an isomorphism in the category of $C^*$-correspondences. Because of that we will  identify $\be(n)$ with $E_n (\La)$
without any further comments.
 The resulting product system of $C^*$-correspondences will be called the product system of the
graph $\La$ and denoted by $\be(\La)$. In what follows we will often view the Dirac functions $\delta_a$ ($a \in \La^0$) and
$\delta_{\la}$ ($\la \in \La^{e_j}$) as elements respectively of $\alg_0$ and of $E_{j}(\La)$.

\subsection*{Representations of $\be(\La)$}

\begin{deft}\label{TCK}
Suppose that $\La$ is a higher-rank graph. A family of partial isometries $\{x_{\la}:\la \in \La\}$
in a $C^*$-algebra $B$ is called a Toeplitz $\La$-family if the following are satisfied:
\begin{rlist}
\item $\{x_a: a \in \La^0\}$ is a family of mutually orthogonal projections;
\item $x_{\la} x_{\mu} = x_{\la \mu}$ if $\la, \mu \in \La$, $s(\la) = r(\mu)$;
\item $x_{\la}^* x_{\la} = x_{s(\la)}$ if $\la \in \La$;
\item if $n \in \Nr\setminus\{0\}$, $a \in \La^0$  and $F\subset\{\la\in\La^n: r(\la)=a\}$ is finite then
$x_a \geq \sum_{\la \in F} x_{\la} x_{\la}^*$;
\end{rlist}
If $\La$ is finitely aligned and additionally the condition \[
 \textrm{(v)} \;\; x_{\mu}^* x _{\nu} = \sum_{\mu
\alpha = \nu \beta \in MCE(\mu, \nu)} x_{\alpha} x_{\beta}^*\hspace*{6.5 cm}\] is satisfied for all $\mu, \nu \in \La$, the
family $\{x_{\la}:\la \in \La\}$ is called a Toeplitz-Cuntz-Krieger family.
\end{deft}

In \cite{toep} isometric representations of $\be(\La)$ are called Toeplitz representations. They are given by Toeplitz families.

\begin{tw}[{\cite{toep}}, Theorem 4.2] \label{isorep}
Let $\La$ be a higher-rank graph. There is a 1-1 correspondence between isometric representations of $\be(\La)$ on a Hilbert
space $\Hil$ and Toeplitz $\La$ families in $B(\Hil)$. The correspondence is given by
\[ \sigma(\delta_{a}) = x_a,\;\;\; a \in \La^0,\]
\[ T^{(j)} (\delta_{\la}) = x_{\la}, \;\;\; j \in \oner, \la \in \La^{e_j}.\]
\end{tw}

Note that $\sigma$ defined as above is nondegenerate if and only if $\sum_{a \in \La_0} x_a =
I_{\Hil}$, where the sum is understood in the strong operator topology.

It is also possible to give an easy characterisation of those isometric representations which are doubly commuting (equivalently,
Nica-covariant).

\begin{lem} [{\cite{toep}}, Proposition 6.4 ] \label{Nicacov}
Let $\La$ be a finitely aligned higher-rank graph. An isometric representation of $\be(\La)$ on a Hilbert space $\Hil$ is doubly
commuting if and only if the corresponding Toeplitz family is a Toeplitz-Cuntz-Krieger family.
\end{lem}

It is easy to see that if $\La$ is finitely aligned then $\be(\La)$ satisfies the normal ordering condition. Note that by Theorem
5.4 of \cite{toep} $\La$ is finitely aligned if and only if $\be(\La)$ is compactly aligned.

We are now ready to define objects which were the main subject of investigation in \cite{graphdil}.

\begin{deft}        \label{Lcont}
Let $\Hil$ be a Hilbert space. A family $\Vcont =\{V_{\lambda}:\la \in \La\}$ of operators in
$B(\Hil)$ is called a $\La$-contraction if the following conditions are satisfied:
\begin{rlist}
\item $\forall_{\la,\mu \in \Lambda,\, s(\la)\neq r(\mu)} \; V_{\la} V_{\mu} =0$;
\item $\forall_{\la,\mu \in \Lambda,\, s(\la) = r(\mu)} \; V_{\la} V_{\mu} = V_{\la \mu}$;
\item $\forall_{n \in \Nr} \; \sum_{\la \in \La^n} V_{\la} V_{\la}^* \leq I$;
\item each $V_a$ ($a \in \La^0$) is an orthogonal projection.
\end{rlist}
All infinite sums here and in what follows are understood in the strong operator topology.
\end{deft}

The definition in \cite{graphdil} was slightly different as  we additionally requested that $\sum_{a\in \La^0} V_a =I$. As
explained in that paper the distinction is not very important: then conditions (ii) and (iii) imply that each $V_a$ for $a \in
\La^0$ is a contractive idempotent, hence a projection (so that in particular (iv) is a consequence of (ii) and (iii)). Further
(i) shows that $V_a V_b = 0$ if $b \in \La^0$ and $a\neq b$. Denoting by $p$ the sum $\sum_{a \in \La^0} V_a$ we see that
$V_{\la} = pV_{\la}p$ (by (i) and (ii)). Therefore even if $\sum_{a\in \La^0} V_a =I$ is not satisfied at the outset, it will be
fulfilled by the obvious $\La$-contraction on $p\Hil$. In condition (iii) above it is enough to assume that the inequalities hold
only for $n$ of the form $e_j$, $j \in \oner$.

We write $V_{\lambda \mu}:=0$ if $s(\la) \neq r(\mu)$. Sometimes we will also write $V_{\emptyset}=
I_{\Hil}$.

The following observation is not very complicated but lies at the heart of this section; it was actually the motivating point for
trying to extend the results of \cite{graphdil} to the framework of representations of product systems of $C^*$-correspondences.

\begin{lem} \label{ccrep}
Let $\La$ be a higher-rank graph. There is a 1-1 correspondence between completely contractive
representations of $\be(\La)$ on a Hilbert space $\Hil$ and $\La$-contractions in $B(\Hil)$. The
correspondence is given by
\begin{equation} \sigma(\delta_{a}) = V_a,\;\;\; a \in
\La^0,\label{rep1}\end{equation}
\begin{equation} T^{(j)} (\delta_{\la}) = V_{\la}, \;\;\; j \in \oner, \la \in
\La^{e_j}.\label{rep2}\end{equation}
\end{lem}

\begin{proof}

Let $\Vcont$ be a $\La$-contraction. Fix for a moment $j \in \oner$ and write $E$ instead of $E_j(\La)$. To show that
$T:=T^{(j)}$ defined by the linear extension of the formula \eqref{rep1} is completely contractive consider  a matrix
$(x_{il})_{i,l=1}^n$ of finitely supported functions in $E$ and let $\xi_1, \ldots, \xi_n$ be vectors in $\Hil$. Let $T_n$ denote
the $n$-th matrix lifting of $T$ and write $\xi = [\xi_1, \cdots, \xi_n]^{\textrm{T}} \in \Hil^{\oplus n}$. Then
\begin{align*} \|T_n \left((x_{i,k})_{i,k=1}^n\right) \xi\|^2 &= \sum_{i,l,k=1}^n \langle T(x_{ik}) \xi_k, T (x_{il})
\xi_l \ra \\&= \sum_{i,k,l=1}^n \left\langle \sum_{\la \in \La^{e_j}} x_{ik}(\la) V_{\la} \xi_k ,
\sum_{\la' \in \La^{e_j}} x_{il}(\la') V_{\la'} \xi_l\right\ra \\&= \sum_{i=1}^n \left\langle
\sum_{\la \in \La^{e_j}} V_{\la} \sum_{k=1}^n  x_{ik}(\la) V_{s(\la)} \xi_k, \sum_{\la' \in
\La^{e_j}} V_{\la'} \sum_{l=1}^n  x_{il}(\la') V_{s(\la')} \xi_l\right\ra.
\end{align*}
Define for each $i =1, \ldots, n$ and $\la \in \La^{e_j}$
\[\zeta^{i}_{\la} = \sum_{k=1}^n x_{ik}(\la) V_{s(\la)} \xi_k.\]
Then
\begin{align*} \|T_n \left((x_{i,k})_{i,k=1}^n\right)\xi\|^2 =  \sum_{i=1}^n \left\langle \sum_{\la \in
\La^{e_j}} V_{\la} \zeta^i_{\la}, \sum_{\la' \in \La^{e_j}} V_{\la'} \zeta^i_{\la'}\right\ra
=\sum_{i=1}^n  \left\| \sum_{\la \in \La^{e_j}}V_{\la} \zeta^i_{\la}\right\|^2.\end{align*} The
condition (iii) in Definition \ref{Lcont} implies that $\| \sum_{\la \in \La^{e_j}}V_{\la}
\zeta^i_{\la}\|^2\leq \sum_{\la \in \La^{e_j}} \|\zeta_{\la}^i\|^2$. Moreover for $\la \in
\La^{e_j}_a$
\[ \|\zeta^{i}_{\la}\|^2 = \|\sum_{k=1}^n x_{ik}(\la) V_a \xi_k\|^2 = \sum_{k,l=1}^n
\ol{x_{ik}(\la)} x_{il}(\la) \langle  V_a \xi_k, V_a \xi_l \rangle\]  so that
\begin{equation} \|T_n \left((x_{i,k})_{i,k=1}^n\right)\xi\|^2 \leq \sum_{i=1}^n \sum_{\la \in \La^{e_j}} \|\zeta_{\la}^i\|^2 =
\sum_{i=1}^n \sum_{a \in \La^0}\sum_{\la \in \La^{e_j}_a}\sum_{k,l=1}^n \ol{x_{ik}(\la)}
x_{il}(\la) \langle  V_a \xi_k, V_a \xi_l \rangle. \label{longlast} \end{equation} Define for each
$a \in \La_0$ a matrix $A_a \in M_n$ by
\[(A_a)_{k,l} = \sum_{i=1}^n \sum_{\la \in \La^{e_j}_a} \ol{x_{ik}(\la)}
x_{il}(\la).\] Note that the  \eqref{opstruct} implies that
\[ \|(x_{i,k})_{i,k=1}^n\|_{M_n(E)} = \sup_{a \in \La^0} \|A_a\|.\]
On the other hand \eqref{longlast} implies that
\[ \|T_n \left((x_{i,k})_{i,k=1}^n\right)  \xi\|^2 \leq \sum_{a \in \La^0} \langle P_a^{(n)} \xi, (A_a \ot I_{\Hil}) P_a \xi\rangle,\]
where $P_a^{(n)}=P_a \oplus \cdots \oplus P_a \in B(\Hil^{\oplus^n})$. Thus finally
\[ \|T_n \left((x_{i,k})_{i,k=1}^n\right)  \xi\|^2 \leq \sup_{a\in \La^0}\{\|A_a\|\} \sum_{a \in \La^0} \| P_a^{(n)} \xi\|^2 =
\sup_{a\in \La^0}\{\|A_a\|\} \|  \xi\|^2,\] so that $T$ extends to a complete contraction from $E$
to $B(\Hil)$.
 The fact that the continuous linear extension
of \eqref{rep1} yields a representation of $\alg_0$ is immediate and then the routine check shows that $(\sigma, T^{(1)}, \ldots,
T^{(r)})$ is a representation of $\be(\La)$.

Conversely, if $(\sigma, T^{(1)}, \ldots, T^{(r)})$ is a representation of $\be(\La)$ we can use formulas \eqref{rep1} and
\eqref{rep2} to define operators $V_{\la}$ for $\la \in \La^0 \cup \bigcup_{j=1}^r \La^{e_j}$. Given any $\mu \in \La$ due to the
factorisation property we can always write it as a concatenation of elements in $ \La^0 \cup \bigcup_{j=1}^r \La^{e_j}$ and
define $V_{\mu}$ as a corresponding composition. The fact that this gives a unique prescription is a consequence of the fact that
$(\sigma, T^{(1)}, \ldots, T^{(r)})$ is a representation of $\be(\La)$; moreover it is easy to check that conditions (i), (ii)
and (iv) of Definition \ref{Lcont} are satisfied (contractive idempotents in $B(\Hil)$ are orthogonal projections). It remains to
check (iii). By (i) and remarks after the definition of a $\La$-contraction it is enough to do it for $n=e_j$ ($j\in \oner$). Let
then $n \in\bn$ and let $\la_1,\cdots, \la_n$ be distinct elements in $\La^{e_j}$. Then the row matrix $[V_{\la_1} \cdots
V_{\la_n}]$ is equal to $T^n((x_{1k})_{k=1}^n)$, where $x_{1k} = \delta_{\la_k}$. It follows easily from \eqref{opstruct} that
$\|(x_{1k})_{k=1}^n\|_{M_n(E_j)} = 1$, and as $T$ is assumed to be a complete contraction we obtain $\|[V_{\la_1} \cdots
V_{\la_n}]\|\leq 1$ and the result follows.
\end{proof}

Note that the representation $\sigma$ of $\alg_0$ associated to a $\La$-contraction $\Vcont$ is
nondegenerate if and only if $\sum_{a \in \La^0} V_a = I_{\Hil}$.

\begin{lem}
Let $\La$ be a rank-r graph and let $\Vcont$ be a $\La$-contraction on a Hilbert space $\Hil$. The representation $\repT$ of
$\be(\La)$ associated with $\Vcont$ is doubly commuting if and only if for all $i,j\in \oner$, $\la \in \La^{e_i}$, $\mu \in
\La^{e_j}$, there is
\[ V_{\la}^* V_{\mu} = \sum_{\alpha \in \La^{e_i}, \beta \in \La^{e_j}, \mu \beta = \la \alpha} V_{\beta} V_{\alpha}^*.\]
\end{lem}

\begin{proof}
Let $i\in \oner$, $\la \in \La^{e_i}$ and $\xi \in \Hil$. Then
\[T^{(i)} (\delta_{\la} \ot \xi) = V_{\la} \xi\]
and it follows that for $\eta \in \Hil$ \begin{equation} (T^{(i)})^* \eta = \sum_{\la \in \La^{e_i}} \delta_{\la} \ot V_{\la}^*
\eta.\label{adjoint} \end{equation}
The equivalence of the conditions in the lemma follows from straightforward computations.
\end{proof}

It would be interesting and nontrivial to analyse how the results of this section extend to topological higher-rank graphs as
discussed for example in \cite{Trent}.

\section{Dilating graph-contractions via dilating representations of associated product systems of
 Hilbert $C^*$-correspondences} \label{dilgraph}

Here we apply the conclusions of the discussions of previous two sections to obtain the dilations of $\La$-contractions to
Toeplitz-type families.

\begin{tw}
Let $\La$ be a rank-2 graph and let $\Vcont$ be a $\La$-contraction on a Hilbert space $\Hil$. There exists a Hilbert space $\Kil
\supset \Hil$ and a $\La$-contraction $\Wcont$ on $\Kil$ consisting of partial isometries forming a Toeplitz family such that for
each $\la \in \La$
\begin{equation} W_{\la}^*|_{\Hil} = V_{\la}^*.\label{dil}\end{equation}
One may assume that $\Kil =
\overline{\textup{Lin}} \{W_{\la} \Hil: \la \in \La\}$. Under this assumption $\sum_{a \in \La^0} W_a = I_{\Kil}$ if $\sum_{a \in
\La^0} V_a = I_{\Hil}$.
\end{tw}

\begin{proof}
Let $\repT$ be the representation of $\be(\La)$ associated with $\Vcont$ by Lemma \ref{ccrep}. From
Theorem \ref{isorep} it follows that any isometric dilation of $\repT$ has to be given by a
Toeplitz family $\Wcont$ such that \eqref{dil} holds. The main statement therefore follows directly
from Theorem \ref{sol1}.
\end{proof}

Before we identify necessary conditions for the representation of $\be(\la)$ associated to a given $\La$-contraction to satisfy
the Brehmer-Solel condition we need to understand the Hilbert spaces involved. Observe that if $\sigma$ is a representation of
$\alg_0$ on a Hilbert space $\Hil$ then for each $n \in \Nr$ the Hilbert space $\be(\la)(n) \ot_{\sigma} \Hil$ is isometrically
isomorphic to the Hilbert space $ \bigoplus_{a \in \La^0} l^2(\La^n_a;P_a \Hil)$, where $P_a= \sigma(\delta_a)$ for each $a \in
\La^0$.

\begin{lem}
Let $r \in \bn$, let $\La$ be a rank-r graph and let $\Vcont$ be a $\La$-contraction on a Hilbert space $\Hil$. Define for each
$u\subset v \subset \oner$ an operator $P_{u,v}$ on the Hilbert space $ \bigoplus_{a \in \La^0} l^2(\La^{e(v)}_a;V_a \Hil)$ via
the continuous linear extension of the formula
\begin{equation} \label{Pdef} P_{u,v} (\delta_{\mu \nu} \ot \xi)  = \bigoplus_{a\in
\La^0} \sum_{\la \in \La^{e(v)}_a} \delta_{\mu\la} \ot V_{\la}^* V_{\nu} \xi,\end{equation}
 where $\mu \in \La^{e(u)-e(v)}, \nu \in
\La^{e(v)}, r(\nu) = s(\mu), \xi \in V_{s(\nu)} \Hil$. Then the associated representation of $\be(\La)$ satisfies the
Brehmer-Solel condition if and only if
\begin{equation}
\forall_{v \subset \oner} \;\;\;\sum_{u\subset v} (-1)^{|u|} P_{u,v} \geq 0.\label{VSolel}\end{equation}
\end{lem}
\begin{proof}
Direct consequence of the remark before the lemma and the formula \eqref{adjoint}.
\end{proof}

\begin{tw}
Let $\La$ be a higher rank graph and let $\Vcont$ be a $\La$-contraction on a Hilbert space $\Hil$.
Suppose that $\Vcont$ satisfies the condition \eqref{VSolel}, where the operators $P_{u,v}$ are
defined by \eqref{Pdef}.  Then there exists a Hilbert space $\Kil \supset \Hil$ and a
$\La$-contraction $\Wcont$ on $\Kil$ consisting of partial isometries forming a Toeplitz family
such that each $W_{\la}^*$ leaves $\Hil$ invariant and for $\la \in \La$, $\mu \in \La$ such that
$\sigma(\la)_j \neq 0$ implies $\sigma(\mu)_j= 0$ ($j \in \oner$)
\[ P_{\Hil}W_{\la}^*W_{\mu}|_{\Hil} = V_{\la}^*V_{\mu}.\]
 One may assume that $\Kil = \overline{\textup{Lin}} \{W_{\la} \Hil: \la
\in \La\}$; under this assumption the family $\Wcont$ is unique up to unitary equivalence.
\end{tw}

\begin{proof}
Let $\repT$ be the representation of $\be(\La)$ associated with $\Vcont$ by Lemma \ref{ccrep}. From Theorem \ref{isorep} it
follows that any isometric dilation of $\repT$ has to be given by a Toeplitz family $\Wcont$ such that \eqref{dil} holds. The
existence of a regular dilation to a Toeplitz family is a consequence of Theorem \ref{sol2}; Lemma \ref{simpreg} implies that
regularity of the dilation can be expressed by a simple formula above.
\end{proof}

The next corollary is a consequence of Theorem 3.15 of \cite{Solk} and Lemma \ref{Nicacov} above.

\begin{cor} \label{correg}
Let $\La$ be a finitely-aligned higher rank graph and let $\Vcont$ be a doubly commuting
$\La$-contraction on a Hilbert space $\Hil$. Then there exists a Hilbert space $\Kil \supset \Hil$
and a $\La$-contraction $\Wcont$ on $\Kil$ consisting of partial isometries forming a
Toeplitz-Cuntz-Krieger family such that each $W_{\la}^*$ leaves $\Hil$ invariant and for $\la \in
\La$, $\mu \in \La$ such that $\sigma(\la)_j \neq 0$ implies $\sigma(\mu)_j= 0$ ($j \in \oner$)
\[ P_{\Hil} W_{\la}^*W_{\mu}|_{\Hil} = V_{\la}^*V_{\mu}.\]
 One may assume that $\Kil = \overline{\textup{Lin}} \{W_{\la} \Hil: \la
\in \La\}$; under this assumption the family $\Wcont$ is unique up to unitary equivalence.
\end{cor}

The following definition was introduced in \cite{graphdil} as a generalisation of the notion of condition `P' suggested in
\cite{PPois}.

\begin{deft} \label{Popcond}
Let $\Vcont$ be a $\La$-contraction and define for $s \in (0,1)$ the defect operator
\begin{equation}\label{Delta} \Pop =
\sum_{\mu \in \La, \, \sigma(\mu) \leq e} (-s^2)^{|\mu|} V_{\mu} V_{\mu}^*.
\end{equation}
The family $\Vcont$ is said to satisfy the Popescu condition (or condition `P') if there exists $\rho\in (0,1)$ such that for all
$s \in (\rho,1)$ the operator $\Pop$ is positive.
\end{deft}

In this context Theorem \ref{main} can be used to establish the following:

\begin{tw}[Theorem 3.1, {\cite{graphdil}}] \label{quoted}
Let $\La$ be a finitely-aligned higher rank graph and let $\Vcont$ be a $\La$-contraction on a Hilbert space $\Hil$ which
satisfies the Popescu condition. Then there exists a Hilbert space $\Kil \supset \Hil$ and a $\La$-contraction $\Wcont$ on $\Kil$
consisting of partial isometries forming a Toeplitz-Cuntz-Krieger family such that $W_{\la}^*|_{\Hil} =V_{\la}^*$ for each $\la
\in \La$.
 One may assume that $\Kil = \overline{\textup{Lin}} \{W_{\la} \Hil: \la
\in \La\}$; under this assumption the family $\Wcont$ is unique up to unitary equivalence.
\end{tw}
\begin{proof}
It follows from the remark stated after Lemma \ref{Nicacov} that $\be(\La)$ has the normal ordering
property. Most of the statements in the theorem follow therefore immediately from Theorem
\ref{main} by now standard applications of the identifications obtained in Section \ref{graphrep}.
The only extra element is double commutativity and uniqueness of the dilation. The first follows
from the fact that if the graph is finitely aligned, then the Fock-Toeplitz representation is
automatically doubly commuting in the natural sense and therefore so is its $^*$-homomorphic image
yielding the dilation in Theorem \ref{main}. The uniqueness follows from the remarks after Theorem
\ref{dc=streg}.
\end{proof}


\begin{thebibliography}{fowl}


\bibitem [Bre] {Brehmer} S.\,Brehmer, \"Uber vetauschbare Kontraktionen des Hilbertschen Raumes, \emph{Acta Sci.\,Math. Szeged} \textbf{22}  (1961), 106--111.

\bibitem[DPY] {DPY} K.\,Davidson, S.C.\,Power and D.\,Yang, Dilation theory for rank 2 graph algebras, \emph{J. Operator
Theory}, to appear.

\bibitem [Fow] {fowl} N.\,Fowler, Discrete product systems of Hilbert bimodules,
  \emph{Pacific J. Math.}  \textbf{204} (2002), 335--375.

\bibitem [FoR] {fowlrae} N.\,Fowler and I.\,Raeburn, The Toeplitz algebra of a Hilbert bimodule,
  \emph{Indiana Univ. Math. J.}  \textbf{48} (1999), 155--181.

\bibitem [GaS] {GaSu} D.\,Ga\c{s}par and N.\,Suciu, On the intertwinings of regular dilations,
  \emph{Ann.\,Pol.\,Math.}  \textbf{LXVI} (1997), 105--121.

\bibitem[KuPa]{kupa}
A.\,Kumjian and D.\,Pask, Higher rank graph $C^*$-algebras, \emph{New York J.\,Math.\,} \textbf{6}
(2000), 1--20.

\bibitem [Lan] {Lance} E.C.\, Lance, ``Hilbert $C^*$-modules", LMS Lecture Notes Series \textbf{210},
Cambridge University Press, Cambridge, 1995.



\bibitem [MS$_1$] {MSgen} P.\,Muhly and B.\,Solel, Tensor algebras over $C^*$-correspondences (Representations, dilations and $C^*$-envelopes),
\emph{J.\,Funct.\,Anal.}  \textbf{158} (1998), 389--457.

 \bibitem [MS$_2$] {MSker} P.\,Muhly and B.\,Solel,
The Poisson Kernel for Hardy Algebras,   \emph{Complex Analysis and Operator Theory}, to appear.

\bibitem[Pop]{PPois}
G.\,Popescu, Poisson transforms on some $C^*$-algebras generated by isometries, \emph{J.
Funct.\,Anal.} \textbf{161} (1999), no.1, 27--61.

\bibitem[PoS]{SteveBaruch}
S.C.\,Power and B.\,Solel, Operator algebras associated with unitary commutation relations, \emph{preprint}, arXiv:0704.0079v1.

\bibitem [Rae] {book} I.\,Raeburn, ``Graph $C^*$-algebras", CBMS Regional Conference Series in Mathematics, 103, Providence, RI, 2005.

\bibitem [RaS] {toep} I.\,Raeburn and A.\,Sims, Product systems of graphs and
the Toeplitz algebras of higher-rank graphs, \emph{J.\,Operator Theory}   \textbf{53}  (2005),  no.
2, 399--429.

\bibitem [RaW] {Morita} I.\,Raeburn and D.\,Williams, ``Morita Equivalence and Continuous-Trace-$C^*$-Algebras", Mathematical Surveys and Monographs, 60. American Mathematical Society, Providence, RI, 1998.

\bibitem[SZ$_1$] {Wold} A.\,Skalski and J.\,Zacharias, Wold decomposition for representations of product systems of $C^*$-correspondences,
 \emph{International Journal of Mathematics} \textbf{19} (2008), no.4, 455-479.

\bibitem[SZ$_2$] {graphdil} A.\,Skalski and J.\,Zacharias, Poisson transform for higher-rank graph algebras and its applications, \emph{J.\,Operator
Theory}, to appear.

\bibitem [So$_1$] {Soltwo} B.\,Solel, Representations of product systems over semigroups and dilations of commuting CP maps,  \emph{J.\,Funct.\,Anal.}  \textbf{235}  (2006),  no. 2, 593--618.

\bibitem [So$_2$] {Solk} B.\,Solel, Regular dilations of representations of product systems, \emph{Math.\,Proc.\,Royal Irish
Soc.}, to appear.

\bibitem [SzF] {Nagy} B.\,Sz.-Nagy and C.\,Foias,
``Harmonic analysis of operators on Hilbert space", North Holland, Amsterdam, 1970.

\bibitem [Tim] {Timot} D.\,Timotin, Regular dilations and models for multicontractions,  \emph{Indiana Univ.\,Math.\,J.}  \textbf{47}  (1998),  no. 2, 593--618.

\bibitem [Yee] {Trent} T.\,Yeend, Topological higher-rank graphs and the $C^*$-algebras of topological 1-graphs,  \emph{in} Operator theory, operator algebras, and
applications,  231--244, \emph{Contemp. Math.} \textbf{414}, Amer. Math. Soc., Providence, RI, 2006.


\end{thebibliography}
\end{document}